\documentclass[12pt,dvips]{article}

\usepackage{amsmath,amsfonts,amssymb,amsthm,graphicx,verbatim,subfig}
\usepackage[all]{xy}

\ifx\pdftexversion\undefined

\usepackage[a4paper,colorlinks,link
=black,filecolor=black,citecolor=black,urlcolor=black,pdfstartview=FitH]{hyperref}
\else

\usepackage[a4paper,colorlinks,linkcolor=black,filecolor=black,citecolor=black,urlcolor=black,pdfstartview=FitH]{hyperref}
\fi

\font\sixbb=msbm6
\font\eightbb=msbm8
\font\twelvebb=msbm10 scaled 1095
\newfam\bbfam
\textfont\bbfam=\twelvebb \scriptfont\bbfam=\eightbb
                           \scriptscriptfont\bbfam=\sixbb
\def\bb{\fam\bbfam\twelvebb}
\newcommand{\Rea}{{\bb R}}

\newcommand{\FF}{{\bb F}}

\newtheorem{theorem}{\bf Theorem}[section]
\newtheorem{claim}[theorem]{\bf Claim}

\newcommand{\enp}{\begin{flushright} $\Box$ \end{flushright}}
\newcommand{\beq}[0]{\begin{equation}}
\newcommand{\enq}[0]{\end{equation}}

\newcommand{\cs}{{\cal S}}

\newcommand{\prob}{{\rm Pr}}

\newcommand{\td}{{\cal T}_d}

\newcommand{\namedref}[2]{\hyperref[#2]{#1~\ref*{#2}}}

\begin{document}

\title{When does the top homology of a random simplicial complex vanish?}
\author{L. Aronshtam\thanks{Department of Computer Science, Hebrew University, Jerusalem 91904,
    Israel. e-mail: liobsar@gmail.com~.} \and  N. Linial\thanks{Department of Computer Science, Hebrew University, Jerusalem 91904,
    Israel. e-mail: nati@cs.huji.ac.il~. Supported by ISF and BSF grants.}
}

\maketitle

\begin{abstract}
Several years ago Linial and Meshulam~\cite{Linial:Meshulam} introduced a model called $X_d(n,p)$ of random $n$-vertex $d$-dimensional simplicial complexes. The following question suggests itself very naturally: What is the threshold probability $p=p(n)$ at which the $d$-dimensional homology of such a random $d$-complex is, almost surely, nonzero? Here we derive an upper bound on this threshold. Computer experiments that we have conducted suggest that this bound may coincide with the actual threshold, but this remains an open question.
\end{abstract}

\section{Introduction}
We study random simplicial complexes in the $X_d(n,p)$ model which was introduced in~\cite{Linial:Meshulam} and further studied in~\cite{MW09, K10, CFK10, BHK11, ALLM}. We quickly recall the basic features of this model. A simplicial complex $X \sim X_d(n,p)$ has $n$ vertices, and a full $(d-1)$-dimensional skeleton. The $d$-dimensional faces of $X$ are selected uniformly and independently with probability $1 \ge p \ge 0$. The parameter $p$ may (and actually will usually) be dependent on $n$. We fix once and for all an arbitrary field $\FF$ and let
$H_i(X)=H_i(X;\FF)$ denote the $i$-th homology group of $X$ with coefficients in $\FF$, and $h_i(X)=\dim_{\FF} H_i(X)$. The main result of this paper is an upper bound on the threshold for the almost sure nonvanishing of $H_d(X)$ for $X \sim X_d(n,p)$. This question was addressed in~\cite{ALLM} and the present paper is, in many ways, a continuation of that paper.

To put this result in perspective, it may be useful to recall the situation for $d=1$ in which case $X_d(n,p)$ coincides with the Erd\H{o}s-R\'enyi $G(n,p)$ model of random graphs. For every $1 > c > 0$, a random graph from $G(n,\frac cn)$ contains a cycle with some probability $1 > f(c) > 0$. Therefore, the threshold for the appearance of a $1$-homology in $G(n, \frac cn)$ is coarse. Namely, the probability that the random graph from $G(n, \frac cn)$ contains a cycle, jumps from a positive number to $1$ as $c$ changes from $1-\epsilon$ to $1$. The behavior for dimension $d > 1$ is similar. This is due to the fact that $\partial \Delta_{d+1}$, the boundary of a $(d+1)$-simplex occurs with positive probability in $X_d(n,\frac cn)$ for every $c>0$. 

However, the situation in one dimension and in higher dimensions are quite different when we consider the size of the first occurring cycle. This issue is a little easier to discuss in the language of a random graph (or complex) process, where edges (resp. $d$-faces) are added at random one at a time. For graphs, the length of the first generated cycle is distributed according to a certain known distribution~\cite{FKP}. In contrast, in the random $d$-dimensional complex process, it follows from~\cite{ALLM} that the first emerging cycle in $H_d(X)$ is either $\partial \Delta_{d+1}$, or it has cardinality $\Omega(n^d)$. We have conducted fairly extensive computer experiments with $d=2$. In these experiments the situation was {\em always} that in the latter case, the first emerging cycle included {\em all} $n$ vertices. Whether this can be proved and whether the same holds for general $d$ remains a subject for further study.

The {\em degree} of a $(d-1)$-face in a $d$-dimensional complex is the number of $d$-faces that contain it. A $(d-1)$-face of degree zero, i.e. one that is contained in no $d$-face is said to be {\it isolated}. A $(d-1)$-face of degree $1$ is said to be {\it free}. The removal of a free $(d-1)$-face and the unique $d$-face that contains it is called an {\it elementary collapse}. We recall~\cite{Hatcher} that an elementary collapse is a homotopy equivalence. Given a complex $X$, we carry out a series of elementary collapses that take place in {\em phases}. At the beginning of a phase we list all $(d-1)$-faces in the complex which are currently free, and we scan them in an arbitrary order. As we arrive at a $(d-1)$-face $\tau$ in the list, one of two things can happen. It may still be free, in which case we apply to it an elementary collapse. It is also possible that when  $\tau$ is reached, it is already isolated, since the unique $d$-face that initially contained it was already eliminated in a previous elementary collapse. In this case we simply skip $\tau$. When we reach the end of the list, the current phase terminates and a new phase commences. 

We denote by $R_i(X)$ the complex gotten from $X$ at the end of phase $i$. In particular $R_0(X)$ is just the randomly drawn complex with which we start. A $d$-face from $R_{i-1}(X) \setminus R_{i}(X)$ is said to be in {\em generation} $i$. A $(d-1)$-face is of {\em generation} $i$ if its degree in $R_{i-1}(X)$ is positive and it either does not belong to $R_{i}(X)$ or is isolated there.

\begin{theorem}
\label{mainTh}
Let $1>\beta=\beta_d>0$ be the unique positive root of the equation
\begin{equation}\label{eq1_thm}
-\ln(1-\beta)=\frac{(d+1)\cdot\beta}{d+1-d\cdot\beta}
\end{equation}
and let $c^{\ast}_d$ be defined as
\begin{equation}\label{eq2_thm}
c^{\ast}_d=\frac{-\ln(1-\beta)}{\beta^d}.
\end{equation}

If $c > c^{\ast}_d$, then a complex $X$ drawn from $X_d(n,\frac{c}{n})$ satisfies asymptotically almost surely
\[
H_d(X) \neq 0.
\]
Specifically, $c_2=2.75381, c_3=3.90708$. Also $\beta_2=0.883414$ and $\beta_3=0.972498$. For large $d$
$$c^{\ast}_d=(d+1)-\frac{d^2+d+1}{\exp(d+1)}+O(\frac{d^2}{\exp(2d)})$$
and $\beta_d=1-\exp(-(d+1))-(1+o_d(1))\frac{(d+1)^2}{\exp(2(d+1))}$.
\end{theorem}

{\bf Note:} Here are a few words about our experiments for $d=2$. We run the random process in which a random $2$-face is added to the complex at each step. The experiment splits according to whether the first cycle to occur is $\partial \Delta_3$ or not. Conditioned on the first cycle not being $\partial \Delta_3$, the numerical estimates that we get for $c_2^{\ast}$ for $n=50, 100, 200$ are  (expectation $\pm$ standard deviation) $2.70424 \pm 0.03115, 2.72886 \pm 0.01534, 2.74149 \pm 0.00733$ respectively. This lends some support to our belief that the bounds attained in Theorem~\ref{mainTh} is the true value of the threshold probability.

The general strategy of our proof is this: An elementary collapse is a homotopy equivalence and in particular it preserves the homology of the complex. Using ideas similar to~\cite{ALLM} we observe what happens as we systematically collapse (in phases, as described above) every free $(d-1)$-face. An elementary collapse eliminates exactly one $(d-1)$-face and one $d$-face. However, it also happens that a non-isolated $(d-1)$-face becomes isolated through collapses on other $(d-1)$-faces. We denote by $C^1_{\ast}$ the (random) set of such faces. Also, let $C^0_{\ast}$ be the set of isolated $(d-1)$ faces in the original $X$. Finally let $C_{\ast} := C^0_{\ast} \cup C^1_{\ast}$ and let $\zeta_{\ast}:=|C_{\ast}|$. Observe that if $f_d(X) > f_{d-1}(X) - \zeta_{\ast}$, then $H_d(X) \neq 0$. Thus the main parts of the proof are these:

\begin{itemize}
\item
Local analysis of the collapsing process.
\item
Computing the expectation $\mathbb{E}(\zeta_{\ast})$.
\item
A measure-concentration argument on the random variable $\zeta_{\ast}$.
\end{itemize}

Our proof uses the fact that every $(d-1)$-face in $C_{\ast}$ corresponds to a zero row in the inclusion matrix of $(d-1)$-faces vs. $d$-faces  (after the collapses). There is another way of establishing the threshold, as done in the upper bound proof in \cite{ALLM}. It is possible to associate to every $(d-1)$-face in $C_{\ast}$ a cocycle in $Z^{d-1}(X)$ and use the Euler-Poincar\'{e} relation to give an upper bound on the threshold. Indeed our proof can be viewed as an extension of the argument of \cite{ALLM}.

As the reader has probably noticed, the above explanation says nothing about the parameter $\beta$ which plays a key role in the theorem. This is done in Section~\ref{overview} below, where we provide a more comprehensive overview of the proof. 

\section{The Probability Space $\td(k,c)$}
\label{s:tree}

We analyze the sequence of $d$-complexes which are obtained, starting from $X$ and repeatedly collapsing, in phases. Our analysis seeks to determine the way at which a given face $\phi$ of dimension $(d-1)$ or $d$ gets collapsed. Note that $\phi$'s generation is completely determined by its local neighborhood in $X$. Concretely, if $\phi$ is of generation $k$, this can be ascertained by observing $\phi$'s radius-$(k+1)$ neighborhood. For every fixed $k$ this neighborhood is almost surely a $d$-tree. We analyze the properties of this neighborhood using an intermediary -- A Galton-Watson-like model of $d$-trees. This model is relatively easy to comprehend, and yet it provides a good approximation to the true local behavior of $X_d(n,\frac cn)$ at the vicinity of $\phi$. This general strategy has been used numerous times, and in particular in \cite{ALLM}.

We turn to provide the necessary definitions. We start with the (recursive) definition of a {\it $d$-tree}. A single $d$-face is a $d$-tree. A $d$-tree on $n+1$ vertices is obtained by taking a $d$-tree $T$ on $n$ vertices and adding to it a new $d$-face $v\cup \tau$ and its $(d-1)$-skeleton. Here $\tau$ is a $(d-1)$-face of $T$, and $v$ is a new vertex.
A {\it rooted $d$-tree} is a $d$-tree in which we designate one $(d-1)$-face to be the root.

Associated with every $d$-complex $Z$ is a graph $G_Z$ whose vertices are
the $d$-faces and the $(d-1)$-faces of $Z$. An edge between a $d$-face and a $(d-1)$-face stands for inclusion, and two $(d-1)$-faces are adjacent when they are contained in a $d$-face of $Z$. We freely apply to $Z$ graph-theoretic notions from $G_Z$ such as distance, diameter and radius from a vertex. Thus the distance between two $(d-1)$-faces in $Z$ is the distance of the two corresponding vertices in $G_Z$.

We define a probability space $\td(k,c)$ of $d$-trees of radius $\le k$ from a $(d-1)$-face $\tau$ that is the root of the tree. Thus $\td(0,c)$ is just the root $(d-1)$-face $\tau$. For $k>0$ we sample a $d$-tree from $\td(k,c)$ as follows:
\begin{itemize}
\item Sample a $d$-tree $T$ from $\td(k-1,c)$.
\item For each $(d-1)$-face $\theta$ in $T$ at distance $k-1$ from $\tau$
\begin{itemize}
\item Sample an integer $j$ from the $\mbox{Poisson}(c)$ distribution.
\item Create $j$ new vertices $t_1,\ldots t_j$ and add $j$ new $d$-faces $\theta\cup t_i$ to $T$ for $i=1,\ldots,j$.
\end{itemize}
\end{itemize}

It is useful now to introduce a variation on the notion of the collapsing process. This is a process which we call {\em $\theta$-collapsing}, where $\theta$ is a $(d-1)$-face. This process is identical to the process of collapsing in phases, except that $\theta$ {\em must not be collapsed}, even when it happens to be free.
We analyze how a random $d$-tree $T\in\td(k+1,c)$ behaves under the $\tau$-collapsing process where $\tau$ is the root of $T$. It is obvious that after $k+1$ phases of $\tau$-collapsing, $T$ collapses to $\tau$, but we need to know whether $\tau$ becomes isolated in phase $k+1$ or sooner. To this end we define the event ${\cal C}_r(k+1,d,c)$ that $\tau$ belongs to generation earlier than $r$, where $r \le k+1$. We denote the probability of ${\cal C}_r(k+1,d,c)$ by $\gamma_r(k+1,d,c)$. As mentioned above, whether or not $\tau$ becomes isolated at time $< r$ depends only on its radius-$r$ neighborhood in $T$. In particular, $\gamma_r(k+1,d,c)=\gamma_r(r,d,c)$ if $k+1\ge r$. We denote $\gamma_r(r,d,c)$ by $\gamma_r(d,c)$. Let us calculate the numbers $\gamma_r(d,c)$ for small $r$. Clearly $\gamma_0(d,c)=0$. Also, ${\cal C}_1(k+1,d,c)$ is the event where $T$ consists only of its root $\tau$, so that
\begin{equation}
\label{gamma}
\gamma_1(d,c)=e^{-c}.
\end{equation}

Notice that $\tau$, the root of a tree becomes isolated before the $r$-th phase iff each $d$-face $\sigma \supset \tau$ satisfies the following condition. There is a $(d-1)$-subface $\sigma \supset\tau'$ that we view as the root of a $d$-tree $T'$ which we $\tau'$-collapse. In the $\tau'$-collapsing of $T'$, the root $\tau'$ becomes isolated before phase $(r-1)$. Let $\pi_j$ be the probability that a $\mbox{Poisson}(c)$ random variable takes the value $j$ we obtain:

\begin{eqnarray}\label{recurgamma}
\gamma_r(d,c)=\sum_{j=0}^{\infty}\pi_j (1-(1-\gamma_{r-1}(d,c))^d)^j=\nonumber\\
\sum_{j=0}^{\infty}\frac{c^j}{j!} e^{-c} (1-(1-\gamma_{r-1}(d,c))^d)^j=\\
\exp(-c(1-\gamma_{r-1}(d,c))^d)~\nonumber.
\end{eqnarray}

We denote by ${\cal B}_k(d,c)$ the event that the root of $T \in \td(k+1,c)$ belongs to a generation later than $k$. The probability of this event is $\beta_k(d,c)$. Clearly

\begin{equation}
\label{defbeta}
\beta_k(d,c)=1-\gamma_{k+1}(d,c).
\end{equation}

\section{The Neighborhood of a $(d-1)$-face}
\label{Neighbor}

The next step is rather standard in arguments of the sort we are making. Most of the necessary details are to be found in~\cite{ALLM}, and we now provide a few additional comments and explanations. The purpose is to show that the Poisson-distribution-based tree considered above approximates arbitrarily closely (as $n \rightarrow \infty$) the actual local behavior of our random complex.

How does the neighborhood of a $(d-1)$-face $\tau$ in a $d$-dimensional complex $X$ look like? The $0$-neighborhood $\cs_0$, consists of $\tau$ alone. The $i$-th neighborhood $\cs_i$ is the complex generated by the $d$-faces in $\cs_{i-1}$, and the additional $d$-faces that contain a $(d-1)$-face in $\cs_{i-1}$. We denote by $v_i$ the number of vertices in $\cs_i$.

Let $A_k$ be the event (in $X_d(n,p)$) that $\cs_k$ is a $d$-tree. Let $D$ be the event that every $(d-1)$-face in $X\sim X_d(n,p)$ has degree $\le \log n$.

The argument in \cite{ALLM} has two parts. One shows first 
\begin{claim}
\label{subtree}
Let $k$ and $c>0$ be fixed and $p=\frac{c}{n}$. Then
$$\prob [A_{k} \cap D]=1-o(1).$$
\end{claim}

The next step is to show that conditioned on the event $A_{k} \cap D$ the following recursive random process generates the typical $k$-neighborhood in $X$. As before $\cs_0=\tau$ and for $i>0$, with $\cs_{i-1}$ already in place, the next layer $\cs_i$ is generated according to the following rule: For each $(d-1)$-face $\theta$ in $\cs_{i-1}$ at distance $i-1$ from $\tau$
\begin{itemize}
\item Sample an integer $j$ from the binomial distribution $B(n-v_{i-1},\frac{c}{n})$.  
\item Create $j$ new vertices $t_1,\ldots t_j$ and add $j$ new $d$-faces $\theta\cup t_i$ to $\cs_i$ for $i=1,\ldots,j$.
\end{itemize} 

The only difference between this random process and the way we defined $\td(k,c)$ is that we sample the integer $j$ from 
$B(n-v_{i-1},\frac{c}{n})$ and not from $\mbox{Poisson}(c)$. Notice that if $X \in A_{k} \cap D$ then $v_k=O(\log^{k}n)$ and the total variation distance between $\cs_{k}$ and $\td(k,c)$ is $o(1)$.

Thus for a $(d-1)$-face $\tau$
\begin{itemize}
\item
$\Pr(\deg_{R_{k}(X)}(\tau)>0) = (1-o(1))\beta_k(d,c)$,
\item
$\Pr(\deg_{R_{k-1}(X)}(\tau)=0) = (1-o(1))\gamma_k(d,c))$.
\end{itemize}

Consider an inclusion $\tau \subset \sigma$ of a $(d-1)$-face and a $d$-face. Let $S'_k$ be the $k$-th neighborhood of $\tau$ in $X\setminus \sigma$. We can apply Claim~\ref{subtree} to $S'_k$ and conclude that with probability $1-o(1)$ it is a $d$-tree in which every $(d-1)$-face has degree at most $\log n$. To randomly generate $S'_k$ we just run the random process that generates $S_k$ and modify it, by excluding $\sigma$ from $S'_1$. Thus

\begin{itemize}
\item
$\Pr(\deg_{R_{k}(X)\setminus\sigma}(\tau)>0) = (1-o(1))\beta_k(d,c)$,
\item
$\Pr(\deg_{R_{k-1}(X)\setminus\sigma}(\tau)=0) = (1-o(1))\gamma_k(d,c))$.
\end{itemize}

\section{An Upper Bound For The Threshold}

As usual, we associate a boundary operator with the $d$-dimensional complex $X$. This linear operator corresponds to an $f_{d-1}(X) \times f_d(X)$ matrix $M$ whose rows and columns are indexed by $X$'s $(d-1)$ resp. $d$-faces. All entries of $M$ are in $\{-1,0,1\}$ and are defined as follows. Every $d$-face $\sigma\in X$ is given an orientation $[v_0,\dots ,v_d]$, and $M_{\sigma\setminus v_i,\sigma}:= (-1)^i$ for every $0\le i\le d$. All other entries of $M$ equal $0$. Since $X$ is $d$-dimensional, $h_d(X)=0$ iff $M$'s right kernel is zero.

Let $N$ be an $a \times b$ matrix and let $\zeta=\zeta(N)$ be the number of zero rows in $N$. Clearly $N$ has a nonzero right kernel if $b > a + \zeta$. We apply this simple observation to $M_i$, the matrix associated with $R_i(X)$. Since an elementary collapse is a homotopy equivalence, $H_d(R_i(X))=H_d(X)$ for all $i$. We conclude that if $s_i(X)\stackrel{\text{def}}{=}f_d(R_i(X))-f_{d-1}(R_i(X))+\zeta_i(X)>0$ for some $i\ge 0$, then $H_d(X)\neq 0$.

Our proof shows that if $c>c_d^{\ast}$ then $X\sim X_d(n,\frac{c}{n})$ satisfies a.s $s(X)=s_{k_{\ast}}(X)>0$. Here $k_{\ast}$ is a large enough constant to be determined later. First we calculate the expectation of $s(X)$, and then show that a.s $s(X)>0$.

\begin{theorem}
\label{rightker}
Let $p=\frac{c}{n}$ with $c> c_d^{\ast}$. Then
$$\prob ~[~X \in X_d(n,p):s(X) > 0]=1-o_n(1)~.$$
\end{theorem}
As the previous discussion shows, this theorem implies Theorem~\ref{mainTh}.

{\bf Proof:} Fix $c$ and $d$ and let $\beta_k,\gamma_k$ stand for $\beta_k(d,c),\gamma_k(d,c)$ resp. and let us fix an arbitrarily small $\epsilon > 0$. Let $\zeta_{\ast}(X)=\zeta (M_{k_{\ast}}(X))$.

A $d$-face $\sigma$ is in $R_{k_{\ast}}(X)$ if $\sigma\in X$ and it is not collapsed in phase $k_{\ast}$ or earlier. In particular, at the end of phase $k_{\ast}-1$, every $(d-1)$-subface of $\sigma$ must be contained as well in a $d$-face other than $\sigma$. Hence 
\begin{equation}
\mathbb{E}[f_d(R_{k_{\ast}}(X))]=(1-o(1)){n\choose d+1}\frac cn\beta_{k_{\ast}-1}^{d+1}
\end{equation}
 
Let $\tau$ be a $(d-1)$-face. Let $\cal{Q}_{\tau}$ be the event that $\tau$ becomes isolated after $k_{\ast}$ $\tau$-collapsing phases. The discussion in Sections \ref{s:tree} and \ref{Neighbor} yields that $\Pr({\cal Q}_{\tau}) = (1-o(1))\gamma_{k_{\ast}+1}$. Let $\cal{P}_{\tau}$ be the event that $\tau$ collapses after $k_{\ast}$ collapsing phases, but does not become isolated after $k_{\ast}$ $\tau$-collapsing phases. Equivalently, $\tau$ becomes a free subface of some $d$-face $\sigma$ before collapsing phase $k_{\ast}$, but all other $(d-1)$-subfaces of $\sigma$ are not free prior to collapsing phase $k_{\ast}$. Consequently,  $\Pr({\cal C}_{\tau}) = (1-o(1))\frac cn n \gamma_{k_{\ast}}\beta_{k_{\ast}-1}^{d}$.

The row corresponding to a $(d-1)$-face $\tau$ becomes a zero row or is removed from the matrix after the collapsing phases if $\tau$ is either isolated or was collapsed. Hence this happens only if the event $\cal{Q}_{\tau}\cup \cal{P}_{\tau}$ occurs. Thus the probability that some $(d-1)$-face belongs to the complex and is not isolated after $k_{\ast}$ collapsing phases is $(1-o(1))(1-(\gamma_{k_{\ast}+1}+c\cdot \gamma_{k_{\ast}}\beta_{k_{\ast}-1}^{d}))$. Therefore 
\begin{equation}
\mathbb{E}[f_{d-1}(R_{k_{\ast}}(X))-\zeta_{k_{\ast}}(X)]=(1-o(1)){n\choose d}(1-(\gamma_{k_{\ast}+1}+c\cdot \gamma_{k_{\ast}}\beta_{k_{\ast}-1}^{d}))
\end{equation}

Consequently,

\begin{eqnarray}
\label{crit_beta}
& &\mathbb{E}[s]= \mathbb{E}[f_d(R_{k_{\ast}}(X))]-\mathbb{E}[f_{d-1}(R_{k_{\ast}}(X))+\zeta_{\ast}]\nonumber\\
& = &(1-o(1)){n\choose d}\left(-\beta_{k_{\ast}}+c\beta_{k_{\ast}-1}^d(1-\beta_{k_{\ast}-1})+c\frac{\beta_{k_{\ast}-1}^{d+1}}{d+1}\right).
\end{eqnarray}

\subsection{overview}
\label{overview}

We are now ready to provide a more detailed explanation of our strategy of proof for Theorem~\ref{mainTh}. We fix an integer $d \ge 2$ once and for all, and some $c > 0$ whose value we discuss below. With this fixed $c$, we obtain a recurrence relation for $\beta_k$. (This follows readily from Equation~(\ref{recurgamma})). Namely, starting with $t=1$ we recurse on $t \rightarrow 1-f(t)$, where $f_c(t)=f(t)=\exp(-c\cdot t^d)$. It is easily verified that for every $c > 0$ the function $1-f(\cdot)$ is increasing in $[0,1]$. As already observed in~\cite{ALLM} there is some $c_{\mbox{collapse}} > 0$ depending only on $d$ so that when $c_{\mbox{collapse}} > c>0$, the only root of $1-f_c(t)=t$ in $[0,1]$ is $t=0$. Therefore, for $c_{\mbox{collapse}} > c>0$ the recurrence $t \rightarrow 1-f(t)$ started at $t=1$ converges to zero.

\begin{figure}[h]
\includegraphics[width=3.3in]{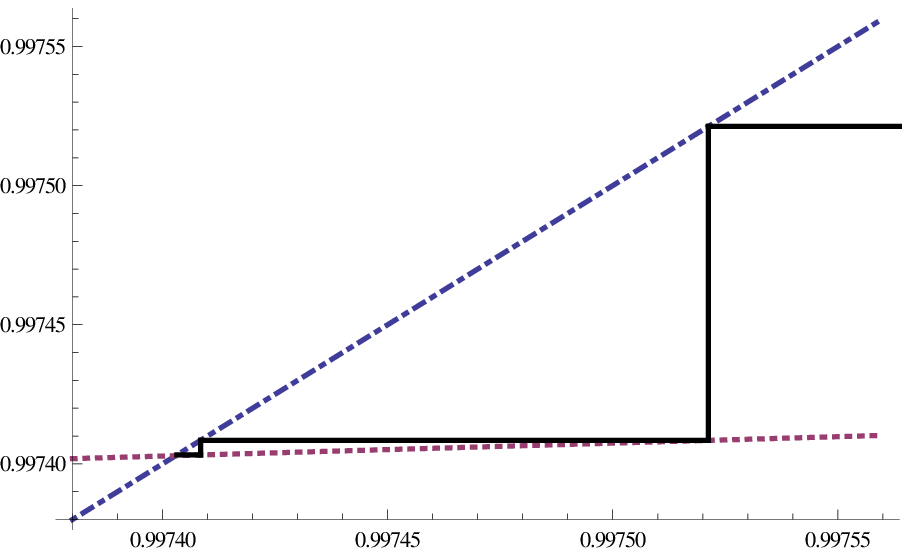}
\includegraphics[width=1.7in]{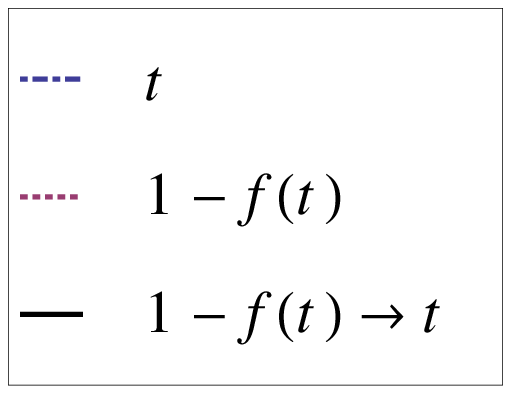}
\caption[Submanifold]{The recurrence relation for $\beta_k$}
\label{recur}
\end{figure}

\begin{figure}[h]
\includegraphics[width=3.2in]{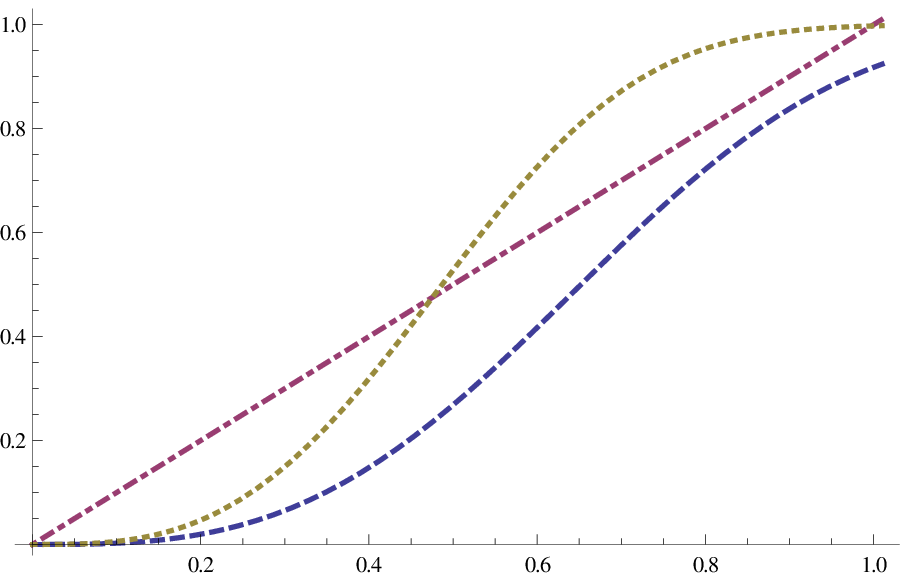}
\includegraphics[width=2.1in]{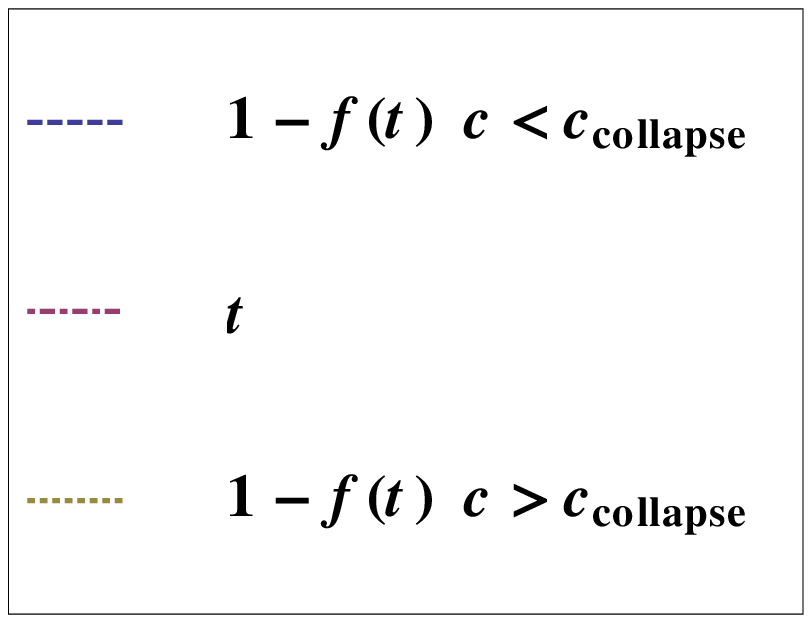}
\caption[Submanifold]{ }
\label{graph}
\end{figure}

Routine calculations (see also Figure~\ref{graph}) yield that for $c > c_{\mbox{collapse}}$ there are exactly two roots in $(0,1)$ to $1-f_c(t)=t$. In this range the above recurrence converges to the larger of these two roots. As Equation (\ref{crit_beta}) shows, $\mathbb{E}(s) > 0$ iff $$-\ln(1-\beta)>\frac{(d+1)\cdot\beta}{d+1-d\cdot\beta}$$ 

In the statement of Theorem~\ref{mainTh} the same calculations are done ``in reverse". Namely, Equation~(\ref{eq2_thm}) states that $1 > \beta > 0$ is a root of $1-f_{c^{\ast}_d}(t)=t$. It only remains to rule out the possibility that Equation~(\ref{eq1_thm}) yields the {\em smaller} root. To this end, note that the solution of Equation~(\ref{eq1_thm}) is $\beta > 1- \exp(-d)$. Moreover, the larger/smaller root of $1-f_c(t)=t$ increases resp. decreases with $c$ and the smaller root is smaller than $1- \exp(-d)$ for every $c > c_{\mbox{collapse}}$.

 Let $k_{\ast}$ be large enough s.t $\beta_{k_{\ast}-1}-\beta_{\ast}<\epsilon$.
Since $c>c_d^{\ast}$ we conclude that for $n$ large enough 
\begin{equation}\label{E(s)}
\mathbb{E}[s] \geq \epsilon' {n\choose d}
\end{equation}
where $\epsilon'>0$ depends only on $c$ and $d$.

\subsection{A concentration of measure argument}

The only missing part of the proof is that $s>0$ almost surely. This is shown using the following version of Azuma's inequality from \cite{McDiarmid}.

\begin{theorem}
\label{mcazuma}
Let $Y_1,\dots,Y_m$ be random variables taking values from $\{0,1\}$. Suppose $\Phi:\{0,1\}^m \rightarrow \Rea$ satisfies:
$$
|\Phi(x)-\Phi(x')| \leq \epsilon_k
$$
for $x$ and $x'$ that differ only at their $k$-th coordinate. Then for every $t>0$:
\begin{equation}
\label{azumeq}
\prob[|\Phi(Y_1,\dots,Y_m)-\mathbb{E}[\Phi(Y_1,\dots,Y_m)]|\ge t]\le2e^{-\frac{2t^2}{\sum_k \epsilon_k^2}}
\end{equation}
\end{theorem}

Let $\sigma_1,\ldots,\sigma_{{n\choose d+1}}$ be the list of all $d$-faces, and let $Y_i$ be the indicator random variable of the event $\sigma_i\in X$. We apply this theorem with $m={{n\choose d+1}}$ and with $\Phi = s$. Let us consider two $d$-complexes $X$ and $X'$ that are identical, except that $\sigma\in X$, but $\sigma\not\in X'$. We need to provide an upper bound on $|s(X)-s(X')|$. In an elementary collapse we eliminate one $d$-face and one $(d-1)$-face. Thus, in particular $f_d(R_i(X))-f_{d-1}(R_i(X))= f_d(X)-f_{d-1}(X)$ for every $i$. Clearly $f_d(X)-f_d(X')=1$ and $f_{d-1}(X)=f_{d-1}(X')$, so we only need a bound on $|\zeta_{\ast}(X')-\zeta_{\ast}(X)|$. As we show $$(d+1)\cdot d^{k_{\ast}} \ge|\zeta_{\ast}(X')-\zeta_{\ast}(X)|.$$

We now compare the collapsing processes as they evolve in $X$ and in $X'$. If a $(d-1)$-face is of generation $i, i'\le k_{\ast}$ in $X, X'$ respectively, then clearly $i \ge i'$. Let $\Theta$ be the set of $(d-1)$-faces in $X$ for which $i > i'$. Clearly $\Theta$ is contained in the $k_{\ast}$-neighborhood of $\sigma$. We classify the faces in $\Theta$ according to their distance from $\sigma$ and show that at distance $j$ from $\sigma$ there are at most $(d+1)\cdot d^{j-1}$ members of $\Theta$. This is clearly true for $j=1$, namely the $d+1$ subfaces of $\sigma$. The general claim is shown by induction on $j$. A $(d-1)$-face $\tau \in \Theta$ at distance $j+1$ from $\sigma$ must have a neighbor (in $G_X$), a $(d-1)$-face $\tau' \in \Theta$ whose distance from $\sigma$ is $j$ so that $\tau'$ is of generation one earlier than $\tau$. In particular, there is a $d$-face $\varphi$ that contains both these $\tau$ and $\tau'$ and is collapsed through $\tau'$.

So let $\tau' \in \Theta$ be a face of generation $\nu'$ in $X'$ whose distance from $\sigma$ is $j$. Since $\tau'$ can collapse only one $d$-face, it can have at most $d$ neighbors in $\Theta$, which are $(d-1)$-faces at distance $j+1$ from $\sigma$. We can conclude that there are at most $(d+1)d^j$ faces in $\Theta$ at distance $j+1$ from $\sigma$.

Clearly $|\zeta_{\ast}(X')-\zeta_{\ast}(X)| \le |\Theta|$, and $|\Theta|\le\sum_{j=1}^{k_{\ast}}(d+1)d^{j-1}\le (d+1)d^{k_{\ast}}$ by the previous discussion. Thus $|s(X)-s(X')|\le (d+1)^{k_{\ast}+1}$, by using (\ref{E(s)}) and (\ref{azumeq}) we conclude:

$$
\prob[s<0]\le2e^{-\frac{2\mathbb{E}^2[s]}{{n\choose d+1} (d+1)^{2k_{\ast}+2}}}=o(1).
$$
{\enp}

\section{Open problems}
\begin{itemize}
\item The most obvious remaining challenge is to determine the correct threshold for the non-vanishing of the $d$-th homology in $X_d(n,p)$. As stated before we believe $c_d^{\ast}$ is this threshold.

\item The present results and those of~\cite{ALLM} strongly suggest that the threshold for collapsibility is substantially smaller than the one for the almost sure non-vanishing of the $d$-th homology. Can one at least show that the two thresholds do not coincide?

\item In the random complex process, what is the distribution of the first emerging cycle in $H_d(X)$? In particular, can one prove that (as suggested by our numerical experiments) it is either $\partial \Delta_{d+1}$ or else it includes all $n$ vertices?

\item
It would be extremely interesting to investigate the inclusion matrices of $(d-1)$-faces vs. $d$-faces of complexes in $X_d(n,\frac cn)$ for values of $c$ between the two thresholds (assuming, of course, that they differ). If the conjectures alluded to in the above questions hold, then this matrix has excellent properties, when viewed as the parity-check matrix of an error-correcting code.

\end{itemize}

\section{Acknowledgment}
This paper could not have been written without our earlier joint work~\cite{ALLM}. We are grateful to Roy and Tomek for their friendship and generosity.

\end{document}